\documentclass[11pt]{article}

\usepackage{amssymb,amsmath,amsthm,rotating}

\newcommand{\n}{\noindent}
\newcommand{\bb}[1]{\mathbb{#1}}
\newcommand{\cl}[1]{\mathcal{#1}}

\theoremstyle{plain}
\newtheorem{lem}{Lemma}

\newtheorem{pro}[lem]{Proposition}

\newtheorem{thm}{Theorem}

\theoremstyle{definition}

\newtheorem{exm}[lem]{Example}

\newtheorem{conj}[lem]{Conjecture}

\begin{document}

\title{A new kind of index theorem\thanks{
Part of this research was done while the author attended the joint India-AMS 
Meeting held in Bangalore in December, 2004. His visit was supported by a grant 
from the DST-NSF Science and Technology Cooperation Programme.}}

\author{Ronald G.\ Douglas\\
\em Dedicated to Krzysztof Wojciechowski on his 50th Birthday}

\date{}
\maketitle

\abstract{Index theory has had profound impact on many branches of mathematics. 
In this note we discuss the context for a new kind of index theorem. We begin, 
however, with some operator-theoretic results. 
In 
\cite{BS} Berger and Shaw established that finitely cyclic 
hyponormal operators have trace-class self-commutators. In \cite{B}, \cite{V} 
Berger and Voiculescu extended this result to operators whose self-commutators 
can be expressed as the sum of a positive and a trace-class operator. In this 
note we show this result can't be extended to operators whose self-commutators 
can be expressed 
as the sum of a positive and a ${\cl S}_p$-class operator.
Then we discuss a conjecture of Arveson 
\cite{Arveson} on homogeneous submodules of the $m$-shift Hilbert space 
$H^2_m$ and propose some refinements of it.

Further, we show how a positive solution would enable one to define $K$-homology 
elements for subvarieties in a 
strongly pseudo-convex domain with 
smooth boundary using submodules of the corresponding Bergman module.  
Finally, we discuss how the Chern character of these classes in cyclic 
cohomology could be defined and indicate what we believe the index to be.}

\setcounter{section}{-1}

\section{Introduction}\label{sec0}

\indent

The complex Hilbert space ${\cl H}$ is said to be a Hilbert module over the 
algebra $A$ if ${\cl H}$ is a unital module over $A$. This is equivalent to a 
representation of $A$ on ${\cl H}$. In the last two decades, there has been 
considerable interest in the study of Hilbert modules for various classes of 
algebras, in part as an approach to multivariate operator theory. In \cite{RP} 
$A$ 
was assumed to be a function algebra and module multiplication to be bounded. 
Other authors (e.g.\ \cite{MS}, \cite{EP}) have considered other kinds of 
algebras. More recently, 
there has been an interest in modules for which $A$ is the algebra of 
polynomials ${\bb 
C}[\pmb{z}]$ with various assumptions such as (1)~~coordinate functions act 
contractively or (2)~~they act as a spherical contraction. In \cite{Arv98} 
Arveson,  
 considered the latter case and identified the $m$-shift space 
$H^2_m$ as having particularly nice properties. In the course of his studies 
\cite{Arveson}, he 
raised a question about the almost reductivity of the submodules of $H^2_m 
\otimes {\bb C}^k$, for $1\le k<\infty$, generated by homogeneous polynomials; 
that is, modules for which the coordinate multipliers and their adjoints have 
compact or $p$-summable cross-commutators. In \cite{A} he established this 
result for 
submodules generated by monomials. (Also, see \cite{arveson4}, \cite{GW}, for 
some subsequent work on this topic.) In \cite{D}, the author extended this 
result 
to a class of commuting weighted shifts which includes the $m$-shift and Bergman 
and Hardy modules for the ball.

In this note we discuss Arveson's conjecture in full generality and  
more. We suggest, in particular,  that submodules of Bergman modules over 
strongly 
pseudo-convex domains of ${\bb C}^m$ with smooth boundary determined by 
subvarieties are $p$-reductive for $p>m$. Moreover, in such a case they 
determine 
odd $K$-homology classes \cite{BDF} for the space equal to the intersection of 
the subvariety with 
the boundary of the domain. Further, one could define a Chern 
character  using the cyclic cohomology of Connes \cite{C}. We conjecture that 
this 
class is  the one determined by the almost complex 
structure 
on the intersection 
of the subvariety and the boundary. Such a result would be a new kind of index 
theorem.

We begin by considering some  results of Berger \cite{B}, which 
extended 
his earlier theorem with Shaw \cite{BS} in operator theory. The latter result 
established 
that 
self-commutators of 
hyponormal operators are trace-class in the presence of finite cyclicity.

My interest in the question of almost reductivity was spurred by Arveson and 
resulted from an 
ongoing dialogue with him on his work on this topic. This rather unusual note 
was the subject of conference talks given in 2005 at IUPUI, Penn State and 
Roskilde University and is presented here to bring to the attention of other 
researchers, what we believe to be most promising and interesting topic.

\section{Results in Operator Theory}\label{sec1}

\indent

Recall that the bounded linear operator $T$ on the Hilbert space ${\cl H}$ is 
said to be {\em hyponormal\/} if the self-commutator $[T^*,T] = T^*T - TT^*$ is 
a positive 
operator. 
In \cite{B} 
Berger and Shaw demonstrated the surprising result that a finitely cyclic 
hyponormal operator has a trace-class self-commutator.  There is also an 
estimate of 
the 
trace involving the degree of cyclicity and the area of the spectrum of $T$ 
but that inequality will not concern us at this time.

Recall that an operator $X$ on ${\cl H}$ is said to belong to ${\cl S}_p$, the 
Schatten--von~Neumann class, for $1\le p<\infty$, if $X$ is compact and the 
eigenvalues of $(X^*X)^{1/2}$ belong to $\ell^p$ (cf.\ \cite{GK}). Now ${\cl 
S}_1$ 
consists of the trace-class operators. One knows that ${\cl S}_p$ is a Banach 
space with dual space ${\cl S}_q$ with $\frac1p +\frac1q = 1$, for $1\le 
p<\infty$, if we identify 
${\cl S}_\infty$ with ${\cl L}({\cl H})$, the space of all bounded operators on 
${\cl H}$. 
Further, ${\cl S}_p$ is a two-sided ideal in ${\cl S}_\infty ={\cl L}({\cl H})$.

In subsequent years,  Berger \cite{B} extended the Berger--Shaw 
Theorem to cover a larger class of operators which is the class we shall 
consider. 
(There was also related work by Voiculescu \cite{V} and Carey--Pincus on this 
class.)
For $1\le p<\infty$, we'll say that an operator $T$ on ${\cl H}$ belongs to 
${\cl A}_p$ if $[T^*,T] = P+C$, where $P\ge 0$ and $C$ is in ${\cl S}_p$. 
Observe that all hyponormal operators 
are 
in 
${\cl A}_p$ as are all operators  $T$ for which $[T^*,T]$ is in ${\cl S}_p$.
Observe also for $p=1$, that there is a well-defined trace on the 
self-commutators 
of 
operators in ${\cl A}_1$ taking 
values in $(-\infty,\infty]$ and that for $T$ in ${\cl A}_1$ we have $[T^*,T]$ 
trace-class iff this trace is finite.

Finally, we will let ${\cl A}_0$ denote the operators $T$ for which 
$[T^*,T]=P+C$ 
with $P$ positive and $C$ compact.

The class of hyponormal operators is  
closed under restriction to invariant subspaces. That is, if $T$ is hyponormal 
and ${\cl V}$ is an invariant 
subspace for $T$, then  $T|_{\cl V}$ is  hyponormal. The 
following lemma shows the same is true for class ${\cl A}_p$. 

\begin{lem}[\cite{B}, \cite{V}]\label{lem1}
If $T$ belongs to ${\cl A}_p$ for $1\le p\le\infty$ or $p=0$, and ${\cl V}$ is 
an 
invariant subspace for $T$, then $T|_{\cl V}$ is in ${\cl A}_p$.
\end{lem}

\begin{proof}
If one writes $T = \left(\begin{smallmatrix} A&B\\ 0&C\end{smallmatrix}\right)$ 
relative to the decomposition ${\cl H} = {\cl V}\oplus {\cl V}^\bot$, and $Q$ 
is 
the orthogonal projection on to ${\cl V}$, then
\[
[(T|_{\cl V})^*, (T|_{\cl V})] =  Q[T^*,T]Q + QTQ^\bot T^*Q,
\]
where $Q^\bot = I-Q$. Since $[T^*,T] =P+C$, with $P\ge 0$ and $C$ in ${\cl 
S}_p$, we have
\[
[(T|_{\cl V})^*, (T|_{\cl V})] = (QPQ + QTQ^\bot T^*Q) + QCQ
\]
and the first sum on the right-hand side is positive while $QCQ$ is in ${\cl 
S}_p$.
\end{proof}

The following result is a special case of a result due to Berger \cite{B}. We 
reproduce the proof since 
it is short and we believe deserves to be better known.

\begin{pro}\label{pro2}
If $T$ is in ${\cl A}_1$ and ${\cl V}$ and $\{{\cl V}_n\}$ are invariant 
subspaces 
for $T$ such that each ${\cl V}_n$ is finite dimensional, ${\cl V}_n\subset 
{\cl 
V}_{n+1}$ for all $n$ and $\bigcup\limits^\infty_{n=1} {\cl V}_n$ is dense in 
${\cl V}$, then $[(T|_{\cl V})^*, (T|_{\cl V})]$ is trace-class.
\end{pro}

\begin{proof}
Let $Q$ and $\{Q_n\}$ be the orthogonal projections onto ${\cl V}$ and $\{{\cl 
V}_n\}$, respectively. Then $\{Q_n\}$ converges in the strong operator topology 
to $Q$. Adopting the same notation as in the preceding proof for the 
representation of the self-commutators, we have $[(T|_{{\cl 
V}_n})^*, (T|_{{\cl V}_n})] = P_n +C_n$ for each $n$ and $[(T|_{\cl V})^*, 
(T|_{\cl 
V})] = P+C$. Moreover, the sequence $\{P_n\}$ converges strongly to $P$ while 
the sequence $\{C_n\}$ converges to $C$ in the norm on ${\cl S}_1$.

Since $T|_{{\cl V}_n}$ is finite rank, we have $Tr[(T|_{{\cl V}_n})^*, 
(T|_{{\cl 
V}_n})] 
= 0$ and hence $0\le TrP_n \le\|C_n\|_1$ for all $n$. Further, we have 
$\|C_n\|_1 \to \|C\|_1$ which implies that $\overline{\rm lim}\ TrP_n  \le 
M<\infty$ 
and 
hence $TrP \le M$ using a variant of Fatou's Lemma. Therefore, $P$ is 
trace-class from which the result 
follows.
\end{proof}

Actually Berger proved a stronger result. Suppose we have another invariant 
subspace ${\cl V}_0$ contained in all the ${\cl V}_n$ so that the dimension of 
${\cl 
V}_n/{{\cl V}_0}$ is finite for all $n$ and $T|_{{\cl V}_0}$ is in ${\cl A}_1$. 
Then 
the preceding argument yields the same conclusion, namely, that the 
self-commutator of  $T|_{\cl V}$ is 
in 
${\cl S}_1$.

We would like to obtain the analogous result for operators in 
$A_p$. In an earlier version of this note we thought we had proved 
it.\footnote{We want to thank W.B.~Arveson for pointing out the mistake in the 
earlier version of this note.}  
Unfortunately, the following example shows it to be false.

\begin{exm}\label{exm3}
Consider the weighted unilateral shift ${S}_n$ defined on $\ell^2$ with the 
standard basis $\{e_k\}^\infty_{k=1}$ so that
\[
{S}_ne_k = \left\{\begin{array}{ll}
\sqrt{\frac{k}n}\ e_{k+1},&1\le k\le n.\\ e_{k+1},&n<k.\end{array}\right.
\]
An easy calculation shows that $\|[{S}^*_n,{S}_n]\|_p = n^{1-p}$ for 
$1\le 
p<\infty$. If ${\cl  V}_n$ is the subspace of $\ell^2$ spanned by 
$\{e_k\}^\infty_{k=n}$, then ${S}_n{\cl V}_n \subset {\cl V}_n$ and 
$\|[({S}_n|_{{\cl V}_n})^*, ({S}_n|_{{\cl V}_n})]\|_p = 1$ for all $n$. 
Moreover, if we set ${S} = \bigoplus\limits^\infty_{n=1} {S}_n$ acting 
on $\bigoplus\limits^\infty_{n=1} \ell^2$ and ${\cl V} = 
\bigoplus\limits^\infty_{n=1} {\cl V}_n$, then $\|[{S}^*,{S}]\|_p < 
\infty$ but $\|[({S}|_{\cl V})^*, ({S}|_{\cl V})]\|_p = \infty$ for all 
$p$, $1<p<\infty$.
\end{exm}

We conclude that ${S}^*$ is in $A_p$, ${\cl V}^\bot_n$ is spanned by 
generalized eigenvectors for ${S}^*$ (hence one can construct the desired 
sequence of finite dimensional approximates for it) but ${S}^*|_{{\cl 
V}^\bot}$ is not $p$ almost reductive. Observe that ${\cl V}^\bot$ and ${\cl V}$ 
are not finitely cyclic for $S^*$ and $S$, respectively.

We now reframe Proposition \ref{pro2} in a setting  which  makes the hypotheses 
more 
transparent.

\setcounter{thm}{-1}

\begin{thm}\label{thm0}
Let $T$ be an operator in ${\cl A}_1$  and 
${\cl V}$ 
be an invariant subspace for $T$ spanned by generalized eigenvectors for 
$T|_{\cl V}$. 
Then $[(T|_{\cl V})^*,(T|_{\cl V})]$ is in ${\cl S}_1$.
\end{thm}

\begin{proof}
Let $\{f_k\}$ be a sequence of generalized eigenvalues for $T|_{\cl V}$ which 
spans ${\cl 
V}$. Further, let ${\cl V}_n$ be the invariant subspace for $T$ generated by 
$\{f_k\}^n_{k=1}$. Then the $\{{\cl V}_n\}$ are nested, finite dimensional and 
their union is dense in ${\cl V}$. The result now follows from 
Proposition \ref{pro2}.
\end{proof}

As we indicated above, our original goal was to extend Proposition \ref{pro2} to 
${\cl A}_p$ and thereby extend  Theorem \ref{thm0} to this class. Unfortunately, 
Example \ref{exm3} 
shows this is impossible.

Since a finite dimensional invariant subspace for $T$ is spanned by the 
generalized 
eigenvectors for it, the hypotheses of the foregoing theorem is the only way to 
fulfill the condition in  Proposition \ref{pro2}. Berger introduced the 
notion of an 
invariant subspace ${\cl V}_0$ being effectually full in ${\cl V}$ by requiring 
 the denseness in ${\cl V}$ of the set of vectors in ${\cl V}$ that some 
nonzero polynomial in $T$ takes into ${\cl V}_0$. This hypothesis enabled him 
to 
satisfy the weaker hypotheses we have mentioned earlier.

A question which presents itself at this point is whether Theorem \ref{thm0} 
might 
hold for all invariant subspaces of $T$ without any restriction. The following 
example shows that this is not the case.

\begin{exm}\label{exm6}
Consider the Bergman space $B^2({\bb D})$ for the unit disk ${\bb D}$ which can 
be defined as the closure of 
the analytic polynomials in $L^2({\bb D})$ relative to planar Lebesgue measure. 
Further, consider the operator $T = M^*_z$, where $M_z$ is multiplication by 
$z$ 
on $B^2({\bb 
D})$. Note that $T$ lies in ${\cl A}_1$.
\end{exm}

Suppose we have invariant subspaces ${\cl M}$ and ${\cl N}$ for $T$ such that 
$0 
\subset {\cl M}\subset {\cl N} \subset B^2({\bb D})$. If  both $T|_{\cl 
M}$ and $T|_{\cl N}$ are in $A_1$,  then the same is true for 
the 
compression of $T$ to the semi-invariant subspace ${\cl N}/{\cl M}$. (This is 
essentially Theorem 1 in \cite{D}.) However, 
one 
knows (cf.\ \cite{ABFP}) that there is a semi-invariant subspace for the Bergman 
shift so
that the restriction of the Bergman shift to it realizes each proper contraction 
operator on a separable Hilbert space. Hence we can obtain a semi-invariant 
subspace  for which the 
self-commutator of the restriction is not even compact.

Here would be a good place to record  a result which is a refinement of Theorem 
1 of \cite{D}.

\begin{thm}\label{thm1}
If ${\cl M}_1$ and ${\cl M}_2$ are essentially reductive modules for the algebra 
$A$ and $X\colon \ {\cl M}_1\to {\cl M}_2$ is a module map having closed range, 
then both $\text{ker } X$ and $\text{ran } X$ are essentially reductive.
\end{thm}

\begin{proof}
If we write ${\cl M}_1 = (\text{ker } X)^\bot \oplus \text{ker } X$ and ${\cl 
M}_2 = \text{ran } X \oplus (\text{ran } X)^\bot$, then $X$ has the form 
$\left(\begin{smallmatrix} X_0&0\\ 0&0\end{smallmatrix}\right)$, where $X_0$ is 
an invertible operator from $(\text{ker } X)^\bot$ onto $\text{ran } X$.

For $T$ any operator between Hilbert spaces, let $\widehat T$ denote its image 
in the corresponding Calkin algebra, that is, modulo the compact operators. 
Since both ${\cl M}_1$ and ${\cl M}_2$ are essentially reductive, for $\varphi$ 
in $A$ the elements $\hat A_\varphi$ and $\widehat B_\varphi$ are normal, where 
$A_\varphi$ and $B_\varphi $ denote the operators defined by module 
multiplication by $\varphi$.

Let $A_\varphi = \left(\begin{smallmatrix} A^{11}_\varphi&0\\ 
A^{21}_\varphi&A^{22}_\varphi\end{smallmatrix}\right)$ and $B_\varphi = 
\left(\begin{smallmatrix} B^{11}_\varphi&B^{12}_\varphi\\ 
0&B^{22}_\varphi\end{smallmatrix}\right)$ be the representations of $A_\varphi$ 
and $B_\varphi$ relative to the decompositions of ${\cl M}_1$ and ${\cl M}_2$, 
respectively. If we consider the images of all these operators in their 
respective Calkin algebras, we can apply the Fuglede--Putnam Theorem to conclude 
that the relationship $\widehat X\hat A_\varphi = \widehat B_\varphi \widehat X$ 
implies that $\hat A_\varphi\widehat X^* = \widehat X^*\widehat B_\varphi$ and 
therefore, we have $\widehat X^*\widehat X \hat A_\varphi = \hat A_\varphi 
\widehat X^* \widehat X$. This equation in turn implies that $\hat 
A^{21}_\varphi \widehat X^*_0 \widehat X_0=0$ and thus $\hat A^{21}_\varphi = 0$ 
since $\widehat X^*_0 \widehat X_0$ is invertible. However, $\hat A^{21}_\varphi 
= 0$ for all $\varphi$ in $A$ means that ${\cl M}_1$ is essentially reductive. 
Working with ${\cl M}_2$ we conclude that $\widehat X\widehat X^* \widehat 
B_\varphi = \widehat B_\varphi \widehat X\widehat X^*$ and hence $\widehat X_0 
\widehat X^*_0 \widehat B^{12}_\varphi = 0$. Again this implies that ${\cl M}_2$ 
is essentially reductive which completes the proof.
\end{proof}

Unfortunately, since no appropriate analogue of the Fuglede--Putnam Theorem is 
known for the $p$-summable case, such a proof won't allow us to conclude 
$p$-reductivity of the kernel and range if ${\cl M}_1$ and ${\cl M}_2$ are. 
However, in \cite{arveson4} Arveson gives such a result for a specific class of 
Hilbert modules and module maps.

\section{Almost Reductive Hilbert Modules}\label{sec2}

\indent

We now show that   Theorem \ref{thm0} enables one to settle a question about 
cross-commutators for some submodules in the multi-variable setting so long as 
they have trace-class 
cross-commutators. Since the conjecture of Arveson 
\cite{Arveson} motivated this study, let us begin by considering it in some 
detail.

Recall that $H^2_m$, the $m$-shift Hilbert space for $1\le m<\infty$, is 
defined 
using the symmetric Fock space and is a module over ${\bb C}[\pmb{z}]$. 
Moreover, 
Arveson showed that multiplication by each coordinate function $Z_i$ acts 
contractively on $H^2_m$ and all cross-commutators 
$[Z^*_i,Z_j]$ lie in ${\cl S}_p$ for $p>m$ and $1\le i,j\le m$.
Arveson conjectured that the restriction operators $Y_i=Z_i|_{\cl S}$ and their 
adjoints also have ${\cl 
S}_p$ cross-commutators for any submodule ${\cl S}$ of $H^2_m\otimes {\bb C}^k$ 
for 
$1\le k<\infty$ generated by 
homogeneous polynomials. Moreover, he established the result for ${\cl S}$ 
generated by monomials. He has also developed \cite{arveson4} a theory of 
``standard Hilbert modules'' in an effort to establish his conjecture. Another 
proof of the  result for monomial submodules was given in 
\cite{D} 
and it also covered certain commuting weighted shifts. Also, Arveson showed that 
the general case for homogeneous submodules of $H^2_m$ for $m=2$ followed from a 
result of Guo \cite{guo}. Finally, a generalization to the case of certain pairs 
of commuting weighted shifts was recently obtained in \cite{GW}.

The simple matrix calculation used in \cite{D} and the proof of Theorem 
\ref{thm1} 
shows that if $T_1$ and $T_2$ are 
two operators on a Hilbert space ${\cl H}$ with ${\cl G}$ an invariant subspace 
for 
both such that $[T_1,T^*_2]$ lies in ${\cl S}_p$ for $1\le p<\infty$ or $p=0$, 
then the compressions $S_i = T_i|_{\cl G}$ have $[S^*_1,S_2]$ in ${\cl 
S}_p$ iff $[R^*_1,R_2]$ is in ${\cl S}_p$ for $R_i = T^*_i|_{{\cl G}^\bot}$. 
Thus we 
can focus on either ${\cl G}$ or ${\cl G}^\bot$.

We consider  the case of commuting weighted shifts 
using the notation of \cite{D}. We have a weight set $\Lambda$ for the index 
set 
$A_m$, $1\le m<\infty$, with the Hilbert space ${\cl M}_\Lambda$ and the 
weighted shifts defined 
by the coordinate functions $Z_i$, $1\le i\le m$. The weight set $\Lambda$ 
satisfies $(*)$ if the shifts are contractive, $(**)$ if all cross-commutators 
of the 
coordinate multipliers and their adjoints are compact, and $(**)_p$ if the 
latter 
operators lie in ${\cl S}_p$.  Actually, in the following result one can replace 
$(*)$ by assuming only $(*)'$ that 
the $Z_i$ are only bounded.

\begin{thm}\label{thm2}
If $\Lambda$ is a weight set satisfying $(*)'$ and $(**)_1$,  ${\cl S}$ is a  
submodule of ${\cl M}_\Lambda \otimes 
{\bb C}_k$, $1\le k<\infty$, so that ${\cl S}^\bot$ is generated by polynomials, 
and $Y_i=Z_i|_{\cl S}$, then the cross-commutators 
$[Y^*_i,Y_j]$ are in ${\cl S}_1$ for $1\le i,j\le m$.
\end{thm}

\begin{proof}
If we set $T = Z^*_i$ for some fixed $i$, then ${\cl S}^\bot$ is invariant for 
$T$. Moreover,  ${\cl S}^\bot$ is spanned by polynomials. Therefore, $T$ and 
${\cl 
S}^\bot$ satisfy the hypotheses of   Theorem \ref{thm0}  which implies that 
$[Y_i,Y^*_i]$ lies in ${\cl S}_1$   for all $1\le i\le m$. Here we 
are using the fact that the self-commutator of the restriction of $Z_i$ to 
${\cl S}$ lies in $S_1$ iff the same is true for the restriction of 
$Z^*_i$ to ${\cl S}^\bot$.

Now if we take $T =Z^*_j + Z^*_k$ for $1\le j \ne k\le m$, then $T$ and ${\cl 
S}^\bot$ again satisfy the hypotheses of Theorem \ref{thm0}. Therefore, we have 
$[Y_j+Y_k, Y^*_j+Y^*_k]$ lies in ${\cl S}_1$. Since $[Y_j,Y^*_j]$ and 
$[Y_k,Y^*_k]$ lie in ${\cl S}_1$, we conclude that the real part of 
$[Y_j,Y^*_k]$ is in ${\cl S}_1$. Repeating the argument for $T = Z^*_j + 
iZ^*_k$, we see that the imaginary part of $[Y_j,Y^*_k]$ is in ${\cl S}_1$ 
which 
completes the proof.
\end{proof}

Let ${\cl P}_n$ denote the subspace of ${\bb C}[z]$ consisting of homogeneous 
polynomials 
of 
degree $n$. If ${\cl S}$ is 
generated by homogeneous polynomials, then $S = \oplus ({\cl S} \cap {\cl 
P}_n)$. This in turn implies that ${\cl S}^\bot = 
\oplus 
({\cl S}^\bot \cap {\cl P}_n)$ and hence ${\cl S}^\bot$ is generated by 
polynomials. Thus  Theorem \ref{thm0} applies to homogeneous submodules. Instead 
of 
assuming that ${\cl S}^\bot$ is generated by polynomials, which are joint 
generalized eigenvectors for the adjoint of coordinate multipliers, we could 
assume more generally that ${\cl S}^\bot$ is spanned by such vectors.

Observe that we can't consider ${\cl M}_\Lambda \otimes \ell^2$, where $\ell^2$ 
is the infinite-dimensional Hilbert space since the cross-commutators on it 
would no 
longer 
be in ${\cl S}_1$. However, if we consider a finite direct sum of block weighted 
shifts 
satisfying 
the analogues of $(*)'$, $(**)_1$, then the result does carry over 
and 
the blocks could be infinite, so long as the cross-commutators are still in 
${\cl 
S}_1$.

While most natural examples of multi-variate Hilbert modules are not 
1-reductive, one can obtain a family of nontrivial examples in the context of 
commuting 
weighted shifts. 

\begin{exm}\label{exm5}
For $m>1$, if the weight set is taken to be:\ 
$\lambda_{\pmb{\alpha}} = \{(1+\alpha_1+\alpha_2 +\cdots+ 
\alpha_m)!\}^{-\delta}$, then ${\cl M}_\Lambda$ is 1-reductive if $\delta 
>\frac{m-1}2$ and the $Z_i$ are in ${\cl S}_2$ if $\delta > \frac{m}2$. Thus 
${\cl M}_\Lambda$ is a nontrivial example of a 1-reductive Hilbert module for 
$\delta$ satisfying $\frac{m}2 \ge \delta > \frac{m-1}2$ and Theorem \ref{thm2} 
applies.
\end{exm}

As we have indicated, originally we had hoped that 
Theorem \ref{thm0} would extend to ${\cl S}_p$, $p>1$, but as Example \ref{exm3} 
indicates, 
this is not the case. Another approach would be to represent either the 
submodule or the corresponding quotient module as the kernel or cokernel of a 
closed module map to which Theorem \ref{thm1} applies. The difficulty here is 
that the module map must have closed range and we know few conditions that 
guarantee that.

 Since most natural examples of multivariate Hilbert 
modules 
are $p$-reductive only for $p>1$, this approach reveals little about the 
validity of Arveson's conjecture in general either for $H^2_m$ or other natural 
examples. Even though that is the case, let us describe what we believe is a 
natural setting for the conjecture.

Let $\Omega$ be a bounded, strongly pseudo-convex domain in ${\bb C}^m$ with 
smooth boundary and 
$B^2(\Omega)$ be the Bergman space, that is, the subset of functions $f$ in 
$L^2(\Omega)$ relative to volume measure for which $\bar\partial f = 0$ taken 
in 
the sense of distributions. One knows \cite{T} that the module action on 
$B^2(\Omega)$ by functions holomorphic on a neighborhood of the closure of 
$\Omega$ is $p$-reductive for $p>m$. That is, cross-commutators of these 
multiplication 
operators and their adjoints lie in ${\cl S}_p$.  For ${\cl Z}$ a variety of 
$\Omega$, 
let 
$B^2_{\cl Z}(\Omega)$ be the functions in $B^2(\Omega)$ that vanish on ${\cl 
Z}$ 
and let ${\cl Q}_{\cl Z}$ be the quotient module $B^2(\Omega)/B^2_{\cl 
Z}(\Omega)$ 
(cf.\ \cite{DM3}). One can show that  ${\cl Q}_{\cl Z}$ is a contractive 
Hilbert 
module 
over $A(\Omega)$ with support in the closure of ${\cl Z}$. Moreover, since
 $B^2(\Omega)$ is a kernel Hilbert space and evaluation at $\pmb{z}$ in 
$\Omega$ is continuous, there is a vector $k_{\pmb{z}}$ in $B^2(\Omega)$ for 
which $f(\pmb{z}) = \langle f,k_{\pmb{z}}\rangle_{B^2(\Omega)}$ for $f$ in 
$B^2(\Omega)$. The vectors $\{k_{\pmb{z}}\}$ are joint eigenvectors for the 
adjoint 
of  the operators defined by
the module action. Moreover, one has that ${\cl Q}_{\cl Z}$ is the closed span 
of 
$\{k_{\pmb{z}}\mid \pmb{z}\in {\cl Z}\}$. Therefore, this example satisfies the 
same kind of hypotheses as in Theorem \ref{thm2}.

More generally, one can see that one could consider any submodule of 
$B^2(\Omega)$ defined as 
the 
orthogonal complement of a collection of eigenvectors $\{k_{\pmb{z}}\}$ and 
their 
partial derivatives, which are also generalized eigenvectors for the adjoint of 
module action. These submodules include the closures of ideals in the
algebra of functions holomorphic on some neighborhood of the closure of 
$\Omega$. 
In particular, one can consider not just the functions that vanish on a 
subvariety but those that vanish to higher order. Moreover, using the result in 
\cite{D} we see that if these submodules are $p$-reductive for $p>m$, then the 
quotient module obtained from them are also $p$-reductive for $p>m$.

Although the evidence for such a result is perhaps scant we are optimistic 
enough to formulate:

\begin{conj}\label{conj6}
If ${\cl S}$ is a submodule of $B^2(\Omega)$ such that ${\cl S}^\bot$ is spanned 
by joint generalized eigenvectors for the adjoint of the operators defined by 
the module action, then both ${\cl 
S}$ 
and ${\cl S}^\bot$ are $p$-reductive for $p>m$.
\end{conj}

This result, even in the multiplicity one case, would be of considerable 
interest. For a submodule obtained as the closure of a principal ideal $I$ in 
${\bb C}[\pmb{z}]$, the result is equivalent to the weighted Bergman space 
defined for the measure $|p|^2 d$ Vol on $\Omega$ being $p$-reductive for $p>m$, 
where $p(\pmb{z})$ is a generator for $I$. However, one might expect, if 
Conjecture \ref{conj6} holds, for the generalization to finite multiplicity to 
also be 
valid.

\begin{conj}\label{conj7}
The same conclusion as in Conjecture \ref{conj6} for submodules of $B^2(\Omega) 
\otimes {\bb C}^k$.
\end{conj}

There is an even stronger result possible which would be very useful in our 
considerations of the following section. (See \cite{DM2} and \cite{DM4} for the 
necessary definitions.)

\begin{conj}\label{conj8}
If ${\cl M}$ is a finite rank, quasi-free, $p$-reductive Hilbert module over 
$A(\Omega)$ and ${\cl S}$ is a submodule for which ${\cl S}^\bot$ is spanned by 
generalized eigenvectors for the adjoint of the operators defined by the module 
action, then ${\cl S}$ and ${\cl 
S}^\bot$ are $p$-reductive.
\end{conj}

It is quite likely that some additional ``regularity'' hypotheses on ${\cl M}$ 
are necessary 
for the last conjecture to hold.

There is another way to frame the final conjecture using a notion introduced in 
\cite{DM2}. 
Recall that a Hilbert module ${\cl M}$ is said to belong to class $(PS)$ if it 
is spanned by the generalized eigenvectors for the adjoint of the operators 
defined by the module action.
\medskip

\n {\bf Conjecture 8$^{\pmb{\prime}}$.}
Let ${\cl H}$ be a finite rank quasi-free, $p$-reductive Hilbert module over the 
algebra $A(\Omega)$. If ${\cl M}$ is a submodule of ${\cl H}$ such that ${\cl 
H}/{\cl M}$ belongs to the class $(PS)$, then ${\cl M}$ is $p$-reductive.

\section{$\pmb{K}$-Homology Classes}\label{sec3}

\indent

Let ${\cl H}$ be a $p$-reductive Hilbert module over the algebra $A$ and ${\cl 
J}({\cl H})$ be the $C^*$-algebra generated by the operators defined by module 
multiplication on ${\cl H}$ and let ${\cl C}({\cl H})$ be the commutator ideal 
in ${\cl J}({\cl 
H})$.
Then ${\cl C}({\cl H})$ consists of compact operators and hence $({\cl J}({\cl 
H}) + {\cl K}({\cl H}))/{\cl K}({\cl H})$ is a commutative $C^*$-algebra. 
Therefore this quotient $C^*$-algebra is isometrically isomorphic to $C(X_{\cl 
H})$ for some 
compact Hausdorff space 
$X_{\cl H}$. In \cite{KD}, it is shown for $A$ a commutative Banach algebra 
that 
$X_{\cl H}$ can be identified with a closed subset of the maximal ideal space 
$M_A$. Similarly, if $A = {\bb C}[\pmb{z}]$ and the module action of the 
coordinate 
functions are all contractive operators, then one can identify $X_{\cl H}$ as a 
closed 
subset of the unit polydisk ${}^{cl}{\bb D}^m$.

In any case, since we have the short exact sequence $0\to {\cl K}({\cl H}) \to 
{\cl J}({\cl H}) + {\cl K}({\cl H})\to C(X_{\cl H})\to 0$,  one always obtains 
an odd $K$-homology element, denoted $[{\cl 
H}]$, 
in $K_1(X_{\cl H})$. While we hope to investigate these classes more thoroughly 
after additional cases of the conjecture have been established, we want to draw 
attention here to a few natural questions and raise a 
few  more
conjectures. Our aim is to show why these are interesting questions. We focus on
the case of  Bergman spaces over strongly 
pseudo-convex domains with smooth boundary.

\begin{thm}\label{thm3}
Let $\Omega$ be a bounded strongly pseudo-convex domain in ${\bb C}^m$ with 
smooth boundary, 
$B^2(\Omega)$ be the Bergman module, ${\cl  Z}$ be a subvariety of $\Omega$, 
$B^2_{\cl Z}(\Omega)$ be the submodule of functions in $B^2(\Omega)$ that vanish 
on 
${\cl Z}$ and ${\cl Q}_{\cl Z}$ be the quotient module $B^2(\Omega)/B^2_{\cl 
Z}(\Omega)$. If ${\cl Q}_{\cl Z}$ is a $p$-reductive module for the algebra 
of 
functions holomorphic on some neighborhood of ${}^{cl}\Omega$, then $[{\cl 
Q}_{\cl 
Z}]$ is in $K_1({\cl Z}\cap \partial\Omega)$.
\end{thm}

\begin{proof}
The only thing requiring proof is the fact that $X_{Q_{\cl Z}} \subseteq {\cl 
Z}\cap 
\partial\Omega$. This follows from the fact that $X_{B^2(\Omega)} = 
\partial\Omega$ and that $B^2_{\cl Z}(\Omega)$ is a Hilbert module over 
$A(\Omega)/A_{\cl Z}(\Omega)$.
\end{proof}

The question arises as to which element of $K_1({\cl Z}\cap \partial \Omega)$ 
is 
obtained. One can show in some cases such as $\Omega = {\bb B}^m$ that it is 
the 
fundamental class, taking multiplicity into account, determined by the complex 
structure on $\Omega$ or the 
spin$^c$-structure on $\partial\Omega$ (or the negative of these classes) 
and I conjecture that this is true in 
general.\footnote{I thank Paul Baum for discussions on how to define such a 
$K$-homology class which is related to our earlier work \cite{BD}.}$^{,}$
\footnote{In 
\cite{GW} the $K$-homology class obtained for homogeneous modules in $B^2({\bb 
B}^2)$ is consistent with this conjecture.} One problem which arises is that 
$\partial\Omega\cap {\cl Z}$ need not 
be a manifold.

One can show by various means that the $K_1$-classes determined by $B^2({\bb 
B}^m)$ and $H^2_m$ are equal. In fact, the same seems to be true for any kernel 
Hilbert module over ${\bb B}^m$ that is essentially reductive. (An argument 
showing this fact would follow from Conjecture \ref{conj8}.) I suspect the 
same thing is true for the $K_1$-classes obtained for a subvariety ${\cl Z}$, 
that is, 
the $K_1$-class doesn't depend on the kernel Hilbert module over $\Omega$ with 
which one 
begins.

Finally, there is one other issue I would like to raise before concluding. We 
will again frame it in the context of submodules of Bergman modules. Although 
one can 
show 
that $B^2_{\cl Z}(\Omega)$ is $p$-reductive for $p>m$, it is not $p$-reductive 
for any smaller $p$. That is, it has the same degree of ``smoothness'' (cf.\ 
\cite{B}) as does $B^2(\Omega)$. However, I don't believe that is  the case 
for 
${\cl Q}_{\cl Z}$. In particular, in \cite{D}, I showed that its smoothness 
depends 
on the dimension of ${\cl Z}$ or the degree of the Hilbert polynomial \cite{DY} 
for 
${\cl Q}_{\cl Z}$. I will formulate one final conjecture, that an analogous 
result holds in general. We state it only for the case of the unit ball.  

\begin{conj}\label{conj9}
Let $I$ be an ideal in ${\bb C}[\pmb{z}]$ and ${\cl S}$ be the submodule 
obtained from 
its closure in $B^2({\bb B}^m)$. Then the quotient module $B^2({\bb B}^m)/{\cl 
S}$ is $q$-reductive for $q>\dim({\cl Z}\cap {\bb B}^m)$, where ${\cl Z}$ is 
the zero 
variety of $I$.\footnote{In \cite{GW} this conjecture is verified in case 
$\dim(Z\cap {\bb B}^m)\le 1$.}
\end{conj}

There is another line of investigation possible here if this conjecture holds. 
If ${\cl Q}_{\cl Z}$ 
is 
$p$-reductive for $q>\dim({\cl Z}\cap {\bb B}^m)$, then it should be possible to 
define a cyclic cohomology 
class 
following Connes \cite{C} which will be the Chern character of $[{\cl Q}_{\cl 
Z}]$. 
One interesting question is how this class  varies when the subvariety 
${\cl Z}$ changes. For example, suppose one considers ${\cl Z}_{c} = 
\{\pmb{z}\in 
\Omega|{p(\pmb{z}}) = c\}$ for $c$ in ${\bb C}$. As one knows, for 
some 
$c,{\cl Z}_{c}\cap \partial\Omega$ will be a manifold while for 
others, it is not. Moreover, there is also the issue of $\partial\Omega\cap 
{\cl 
Z}_{c}$ being a manifold while ${\cl Z}_{c}$ has singularities in 
$\Omega$.

One fascinating example to consider would be the presentation of the exotic 
spheres found by Brieskorn \cite{brieskorn}. Recall that he exhibited 
analytic polynomials for which an exotic sphere is obtained from the 
intersection of the zero variety of the polynomial in ${\bb C}^n$ with spheres 
of small diameter. Although the precise polynomials he used are not homogeneous, 
this example indicates that one is likely to obtain interesting varieties in our 
context.

 I 
believe  different techniques will be needed to establish such a conjecture. The 
result in 
\cite{DV} provides a lower bound on $p$ if $[{\cl Q}_{\cl Z}]$ is indeed a 
fundamental class for $\partial\Omega\cap {\cl Z}$.

\end{document}